\documentclass[reqno,pdftex,11pt]{amsart}

\usepackage[top=25truemm,bottom=25truemm,left=25truemm,right=25truemm]{geometry}

\usepackage{amsmath,amsthm,amssymb,amsfonts,amscd,mathtools,color,tikz,tikz-cd}
\usepackage{graphicx} 
\usetikzlibrary{arrows.meta}
\numberwithin{equation}{section}

\theoremstyle{plain} 
	\newtheorem{thm}{Theorem}[section]
	\newtheorem*{thm*}{Theorem}
	
	\newtheorem{lem}[thm]{Lemma}
	
	\newtheorem{prop}[thm]{Proposition}
	
	\newtheorem*{conj*}{Conjecture}
	
\theoremstyle{definition}
	\newtheorem{defn}[thm]{Definition}

\theoremstyle{remark}
	\newtheorem{rem}[thm]{Remark}
	
	\newtheorem*{pf}{Proof}
\def\CC{{\mathbb C}}

\def\PP{{\mathbb P}}

\def\ZZ{{\mathbb Z}}

\def\D{{\mathcal D}}
\def\E{{\mathcal E}}

\def\L{{\mathcal L}}

\def\R{{\mathcal R}}
\def\S{{\mathcal S}}
\def\T{{\mathcal T}}

\def\Aut{{\rm Aut}}

\def\End{{\rm End}}

\def\FEC{{\rm FEC}}

\def\Hom{{\rm Hom}}

\def\RHom{{\rm \mathbf{R}Hom}}

\def\mod{{\rm mod}}

\begin{document}
\title[Non-transitivity of the braid group action]{A simple algebraic proof of the non-transitivity of the braid group action on full exceptional sequences}
\date{\today}
\author{Atsuki Nakago} 
\address{Department of Mathematics, Graduate School of Science, Osaka University, 
Toyonaka Osaka, 560-0043, Japan}
\email{u400778f@ecs.osaka-u.ac.jp}
\author{Atsushi Takahashi}
\address{Department of Mathematics, Graduate School of Science, Osaka University, 
Toyonaka Osaka, 560-0043, Japan}
\email{takahashi@math.sci.osaka-u.ac.jp}

\begin{abstract}
Recently, Chang--Haiden--Schroll shows that the braid group action on full exceptional collections in a triangulated category
is not transitive but has infinitely many orbits in general. 
Their proof is based on a geometric model and the theory of branched coverings such as Birman--Hilden theory. 
This paper provides a simple algebraic proof of their theorem.
\end{abstract}
\maketitle
\section{Introduction}
Exceptional collections were first introduced by Gorodentsev--Rudakov~\cite{GR} in the study of coherent sheaves on the projective space $\PP^n$.
Building on their geometric constructions, Bondal~\cite{Bo} developed a more general abstract theory of exceptional objects and mutations for arbitrary triangulated categories.
A \emph{full exceptional collection} provides a semi-orthogonal decomposition of the triangulated category into simple building blocks,
and hence plays a central role in the study of derived categories.
In particular, the existence and classification of such collections have been fundamental in the development of homological mirror symmetry and representation theory.
In singularity theory, distinguished bases of vanishing cycles play a central role. 
In particular, a fundamental subject of study is the derived category of the directed Fukaya category associated to a singularity, and a chosen distinguished basis of vanishing cycles gives rise to a full exceptional collection in this category (see \cite{Se1, Se2}). 
Consequently, the essential importance of working with full exceptional collections has been recognized since the early days of singularity theory.

The braid group acts on the set of full exceptional collections via mutations, and Bondal--Polishchuk~\cite{BP} conjectured that this action, together with the 
action of the translations on objects, is transitive. 
This conjecture has been confirmed in several important cases, including hereditary algebras~\cite{CB, Ri}, the projective plane~\cite{Ru}, weighted projective lines~\cite{Me}, all del Pezzo surfaces~\cite{KO}, and the Hirzebruch surface of degree 2~\cite{IOU}.
However, recently Chang--Haiden--Schroll~\cite{CHS} showed that transitivity does not hold in general:
\begin{thm}[{\cite[Theorem~5.4]{CHS}}]\label{thm: main}
Let $\mu\ge 4$ and $(Q, I)$ the quiver with relations given by
\begin{equation*}
Q:
\begin{tikzcd}
\underset{1}{\bullet} \ar[r, shift left, "{\alpha_1}"] \ar[r, shift right, "{\beta_1}"'] 
& \underset{2}{\bullet} \ar[r, shift left, "{\alpha_2}"] \ar[r, shift right, "{\beta_2}"'] 
& \cdots \ar[r, shift left, "{\alpha_{\mu-1}}"] \ar[r, shift right, "{\beta_{\mu-1}}"'] 
& \underset{\mu}{\bullet}
\end{tikzcd},\quad
I = \langle \alpha_{i}\beta_{i+1},\ \beta_{i}\alpha_{i+1}\ |\ i=1,\dots,\mu-2 \rangle.
\end{equation*}
Let $\CC Q/I$ be the path algebra of $(Q,I)$ and $\D^{b}{\rm mod} (\CC Q/I)$ the bounded derived category of finite dimensional $\CC Q/I$-modules.
The $B_\mu\ltimes \ZZ^\mu$-action on the set of full exceptional collections in $\D^b\mod(\CC Q/I)$ is not transitive but has infinitely many orbits.
\end{thm}
They prove this theorem based on a geometric model and the theory of branched coverings such as Birman--Hilden theory. 
Our goal in this paper is to provide an algebraic proof of their theorem without using a geometric model.

\bigskip
\noindent
{\bf Acknowledgements.}
A.~T. is supported by JSPS KAKENHI Grant Number JP21H04994.

\section{Preliminaries}\label{sec: preliminaries}
This section prepares the terminology and some properties that will be needed later.

\subsection{Full exceptional collections}
Let $\T$ be a triangulated $\CC$-linear category with the translation functor $[1]$. 
For $p\in\ZZ$ the $p$-times composition of $[1]$ will be denoted by $[p]$.
We assume that $\T$ is of {\it finite type}, namely, we have
\[
\sum_{p\in\ZZ}\dim_\CC \Hom_\T(X,Y[p]) <\infty,\quad X,Y\in\T
\]
and that $\T$ is the derived category of an abelian $\CC$-linear category.
Let $\D^{b}{\rm mod}(\CC)$ be the bounded derived category of finite dimensional $\CC$-modules. 
Since $\T$ has a canonical differential graded (dg) enhancement, we have the functor
\[
\RHom_{\T}(-,-):\T^{op}\times \T\longrightarrow \D^{b}{\rm mod}(\CC),
\]
with the property that $H^p(\RHom_{\T}(X,Y))=\Hom_\T(X,Y[p])$ for all $X,Y\in\T$, $p\in\ZZ$.
We also have functorial mapping cones for all morphisms in $\T$.
 
\begin{defn}
Let the notations be as above.
\begin{enumerate}
\item
An object $E \in \T$ is called {\it indecomposable} if $E$ is nonzero and $E \cong X \oplus Y$ 
for some $X,Y\in\T$ then $X$ or $Y$ is the zero object.
\item
An object $E \in \T$ is called {\it exceptional} if $\RHom_\T(E,E) =\CC\cdot {\rm id}_E$ in $\D^{b}{\rm mod}(\CC)$.
If an object is exceptional then it is indecomposable.
\item 
An ordered collection $\E = (E_1, \dots, E_\mu)$ of exceptional objects $E_1, \dots, E_\mu\in\T$ is called 
an {\it exceptional collection} if $\RHom_\T(E_i, E_j) \cong  0$ in $ \D^{b}{\rm mod}(\CC)$ if  $i > j$.
\item 
An exceptional collection $\E$ is called {\it full} if the smallest full triangulated subcategory 
$\langle E_1,\dots, E_\mu\rangle$ of $\T$ containing $E_1,\dots, E_\mu$ is equivalent to $\T$.
\item 
Two full exceptional collections $\E=(E_1, \dots, E_\mu)$, $\E'=(E'_1, \dots, E'_\mu)$ in $\T$ 
are said to be {\it isomorphic} if $E_i \cong E'_i$ for all $i = 1, \dots, \mu$. 
Denote by $\FEC(\T)$ the set of isomorphism classes of full exceptional collections in $\T$
\item
A full exceptional collection $\E= (E_1, \dots, E_\mu)$ is called {\it strong} 
if $\RHom_\T(E_i, E_j)$ is isomorphic in $ \D^{b}{\rm mod}(\CC)$ to a complex concentrated in degree zero,
that is, $\Hom_\T(E_i,E_j[p])=0$ for all $p\ne 0$.
\end{enumerate}
\end{defn}

\begin{defn}
Let the notations be as above.
\begin{enumerate}
\item 
An auto-equivalence $\Phi$ of $\T$ is an exact functor $\Phi:\T\longrightarrow\T$ 
which is also an equivalence of $\T$.
Denote by $\Aut(\T)$ the group of isomorphism classes of auto-equivalences of $\T$.
\item A {\it Serre functor} $\S_\T$ is an auto-equivalence of $\T$ with bi-functorial isomorphisms $\Hom_\T(X,Y) \cong\Hom_\T(Y,\S_\T(X))^*$ 
for all $X,Y\in\T$, where $(-)^*$ denotes $\Hom_{\mod(\CC)}(-,\CC)$.
\end{enumerate}
\end{defn}
The notion of Serre functors is introduced by Bondal--Kapranov~\cite[Section~3]{BK}. 
If a Serre functor exists, it is unique up to canonical isomorphism which commutes with the above bi-functorial isomorphisms.
Therefore, it is usually referred to as {\it the} Serre functor of $\T$ and is denoted by $\S_\T$.
Several properties of $\S_\T$ are given there, among which the following are particularly used in this paper.
\begin{prop}[{\cite[Section~3]{BK}}]
Let the notations be as above.
Assume that $\T$ admits a full exceptional collection $(E_1,\dots,E_\mu)$. 
\begin{enumerate}
\item 
The Serre functor $\S_\T$ exists. 
\item
The Serre functor $\S_\T$ commutes with any auto-equivalence of $\T$, namely, $\S_\T\circ \Phi\cong \Phi\circ \S_\T$ for any $\Phi\in \Aut(\T)$. 
\end{enumerate}
\qed\end{prop}
In this paper, we consider as $\T$ the bounded derived category of finite dimensional modules over a finite dimensional $\CC$-algebra, 
so the first statement is a consequence of \cite[Chapter~1, Section~4.6]{Ha}.

\subsection{Group actions on the set of full exceptional collections}

The Artin {\it braid group} $B_\mu$ on $\mu$-strands is a group presented by the following generators and relations: 
\begin{equation*}
B_{\mu}:=\langle {\boldsymbol\sigma}_1,\dots, {\boldsymbol\sigma}_{\mu-1}\,|\,{\boldsymbol\sigma}_{i}{\boldsymbol\sigma}_{i+1}{\boldsymbol\sigma}_{i}={\boldsymbol\sigma}_{i+1}{\boldsymbol\sigma}_{i}{\boldsymbol\sigma}_{i+1}\text{ for }i=1,\dots, \mu-2,\ 
{\boldsymbol\sigma}_{i}{\boldsymbol\sigma}_{j}={\boldsymbol\sigma}_{j}{\boldsymbol\sigma}_{i}\text{ for }|i-j|>1\rangle.
\end{equation*}
Consider the semi-direct product $B_\mu\ltimes \ZZ^{\mu}$ of $B_\mu$ and the free abelian group $\ZZ^\mu$ of rank $\mu$
where the homomorphism $B_\mu\longrightarrow \Aut(\ZZ^\mu)$ is the composition of a group homomorphism
$B_\mu\longrightarrow {\mathfrak{S}}_\mu$ into the symmetric group ${\mathfrak{S}}_\mu$ and the natural action of ${\mathfrak{S}}_\mu$ on $\ZZ^\mu$.

The group $B_\mu\ltimes \ZZ^{\mu}$ acts by mutations and translations on the set $\FEC(\T)$ of full exceptinal collections in $\T$ from the right~\cite{BP}:
\[
(E_{1}, \dots, E_{\mu})\cdot {\boldsymbol\sigma}_{i}:= (E_{1}, \dots, E_{i-1}, \L_{E_{i}} E_{i+1}, E_{i}, E_{i+2},\dots, E_{\mu}),
\]
\[
(E_{1}, \dots, E_{\mu})\cdot {\boldsymbol\sigma}^{-1}_{i} := (E_{1}, \dots, E_{i-1}, E_{i+1}, \R_{E_{i+1}} E_i, E_{i+2}, \dots, E_{\mu}),
\]
\[
(E_{1}, \dots, E_{\mu})\cdot {\bf n} := (E_{1}[n_1],\dots, E_{\mu}[n_\mu]),\quad {\bf n}:=(n_1,\dots,n_\mu)\in \ZZ^{\mu}.
\]
Here, $\L_{E_i}E_{i+1},~\R_{E_{i+1}}E_{i}$ for $i=1,\dots,\mu-1$ are defined by the following exact triangles
\[
\L_{E_{i}}E_{i+1}\longrightarrow \RHom_\T(E_{i},E_{i+1})\otimes E_{i}\longrightarrow E_{i+1}\longrightarrow \L_{E_{i}}E_{i+1}[1],
\]
\[
\R_{E_{i+1}}E_{i}[-1]\longrightarrow E_{i}\longrightarrow \RHom_\T(E_{i},E_{i+1})^*\otimes E_{i+1}\longrightarrow \R_{E_{i+1}}E_{i}.
\]
\begin{rem}
We consider the right $B_\mu$-action since it is natural to regard this $B_\mu$ as the fundamental group of a configuration space 
of $\mu$ points in the complex plane $\CC$.
\end{rem}

On the other  hand, the group $\Aut(\T)$ of auto-equivalences on $\T$ acts on $\FEC(\T)$ from the left by
\[
\Phi \cdot \E := (\Phi(E_1),\dots,\Phi(E_\mu)), \quad \Phi \in \Aut(\T),\ \E=(E_1,\dots,E_\mu)\in\FEC(\T).
\]
Note that on $\FEC(\T)$ the left action of $\Aut(\T)$ and the right action of $B_\mu\ltimes \ZZ^\mu$ commute, namely, 
\[
(\Phi\cdot\E)\cdot ({\boldsymbol\sigma},{\bf n})=\Phi\cdot(\E\cdot ({\boldsymbol\sigma},{\bf n})), \quad \Phi\in\Aut(\T), \ ({\boldsymbol\sigma},{\bf n})\in 
B_\mu\ltimes \ZZ^\mu.
\]
Of particular importance in this paper are the auto-equivalences given based on spherical objects.
\begin{defn}[{\cite[Definition~2.14]{ST}}]
An object $S\in \T$ is \textit{$m$-spherical} if the following holds
\begin{enumerate}
\item $\RHom_\T(S,S)\cong\CC\oplus\CC[-m]$,
\item For the Serre functor $\S_\T$ of $\T$, $\S_\T(S)\cong S[m]$.
\end{enumerate}
\end{defn}
\begin{prop}[{\cite[Theorem~2.10]{ST}}]
Let $S\in\T$ be an $m$-spherical object. There exists an auto-equivalence $T_{S}$ on $\T$ given by the following exact triangle
\[
T_{S}(X)[-1]\longrightarrow \RHom_\T(S,X)\otimes S\longrightarrow X\longrightarrow T_{S}(X),\quad X\in\T.
\]
\qed\end{prop}

\subsection{}
Let $\mu$ be a positive integer and $(Q, I)$ the quiver with relations given by
\begin{equation*}
Q:
\begin{tikzcd}
\underset{1}{\bullet} \ar[r, shift left, "{\alpha_1}"] \ar[r, shift right, "{\beta_1}"'] 
& \underset{2}{\bullet} \ar[r, shift left, "{\alpha_2}"] \ar[r, shift right, "{\beta_2}"'] 
& \cdots \ar[r, shift left, "{\alpha_{\mu-1}}"] \ar[r, shift right, "{\beta_{\mu-1}}"'] 
& \underset{\mu}{\bullet}
\end{tikzcd},\quad
I = \langle \alpha_{i}\beta_{i+1},\ \beta_{i}\alpha_{i+1}\ |\ i=1,\dots,\mu-2 \rangle.
\end{equation*}
Let $\CC Q/I$ be the path algebra of $(Q,I)$, which is a finite dimensional $\CC$-algebra of global dimension $\mu-1$ (see \cite[Section~4]{D}).
Denote by ${\rm mod}(\CC Q/I)$ the abelian category of finite dimensional $\CC Q/I$-modules and 
by $\D^{b}{\rm mod} (\CC Q/I)$ the bounded derived category of ${\rm mod}(\CC Q/I)$.
In order to simplify the notation, set $\D:=\D^{b}{\rm mod} (\CC Q/I)$.  
The Serre functor on $\D$ is denoted by $\S_\D$.

For each $i=1,\dots, \mu$, denote by $P(i)$ and $I(i)$ the indecomposable projective and injective $\CC Q/I$-modules at the vertex $i$ in $Q$, respectively.
The following well-known properties will be used later.
\begin{prop}
Let the notations be as above.
\begin{enumerate}
\item 
For each $i,j=1,\dots, \mu$, $\Hom_\D(P(i),P(j))$ has a basis consisting of all paths from the vertex $j$ to the vertex $i$ in $Q$.
\item
For all $i=1,\dots, \mu$, the $\CC Q/I$-module $P(i)$ is an exceptional object in $\D$.
\item
The ordered collection $\E_P:=(P(\mu),\dots,P(1))$ of exceptional objects is a full strong exceptional collection in $\D$.
\item
For all $i=1,\dots, \mu$, $\S_\D(P(i))\cong \Hom_\D(P(i),\CC Q/I)^* \cong  I(i)$.
\end{enumerate}
\qed\end{prop}

\begin{prop}\label{prop:1-sph}
Let the notations be as above.
For each $i=1,\dots, \mu-1$, an object $S_{i}\in\D$ defined by the following exact triangle
\begin{equation}\label{1-sph}
P(\mu-i+1)\stackrel{f_{i}}{\longrightarrow} P(\mu-i)\longrightarrow S_{i}\longrightarrow P(\mu-i+1)[1],\quad f_{i}:=\alpha_{\mu-i}+\beta_{\mu-i},
\end{equation}
is $1$-spherical and we have 
\begin{equation}
\RHom_\D(S_{i},S_{j})\cong 
\begin{cases}
\CC\oplus \CC[1] & j-i=0,\\
\CC & j-i=1,\\
\CC[-1] & j-i=-1,\\
0& |i-j|>1.
\end{cases}
\end{equation}
\end{prop}
\begin{pf}
It follows from a direct calculation.
\qed\end{pf}
\begin{prop}\label{prop:l-r braid}
Let the notations be as above. In $\FEC(\D)$, we have 
\[
T_{S_{i}}^{-1}\cdot\E_P =\E_P\cdot{\boldsymbol\sigma}_i,\quad i=1,\dots, \mu-1.
\]
\end{prop}
\begin{pf}
It suffices to check the case $i=1$, since the case of general $i$ is obtained by the same argument 
with $(P(\mu),P(\mu-1))$ replaced by $(P(\mu-i+1),P(\mu-i))$.
Applying $\RHom_\D(-,P(k))$ to \eqref{1-sph} for $i=1$ we obtain
\[
\RHom_\D(S_{1},P(k))\cong
\begin{cases}
\CC[-1] & k=\mu,\mu-1,\\
0 & k\ne \mu,\mu-1.
\end{cases}
\]
Hence $T_{S_{1}}(P(\mu))\cong P(\mu-1)$ and $T_{S_{1}}(P(\mu-2))\cong P(\mu-2)$.
Moreover, we have $\RHom_\D(T_{S_{1}}(P(\mu-1)),\R_{P(\mu-1)}P(\mu))\cong\CC$.
Since all other terms in the full exceptional collection coincide, 
this implies $T_{S_{1}}(P(\mu-1))\cong\R_{P(\mu-1)}P(\mu)$.
Therefore
\[
T_{S_{1}}\cdot\E=(P(\mu-1),\R_{P(\mu-1)}P(\mu),P(\mu-2),\dots,P(1))
=\E_P\cdot\boldsymbol\sigma_1^{-1}.
\]
\qed\end{pf}
Proposition~\ref{prop:1-sph} implies that we have the left $B_\mu$-action on the set $\FEC(\D)$ via
the group homomorphism 
\[
B_\mu\longrightarrow \Aut(\D),\quad {\bf \boldsymbol\sigma}_i\mapsto T_{S_{i}}^{-1},
\]
by Seidel-Thomas \cite[Theorem~2.17]{ST}, and Proposition~\ref{prop:l-r braid} gives an identification of this left $B_\mu$-action with the right $B_\mu$-action on $\FEC(\D)$.

\section{Algebraic proof of the non-transitivity}

In this section, we give an algebraic proof of Theorem~\ref{thm: main} which is proven by 
Chang--Haiden--Schroll~\cite[Theorem~5.4]{CHS} by a geometric method.
We will use the terminology prepared in Section~2.
\subsection{The case $\mu=4$}

We first prove Theorem~\ref{thm: main} for the case $\mu=4$.

\begin{lem}\label{lem: S_+-}
Let $S_+$ and $S_-$ be $\CC Q/I$-modules defined by
\[
S_+=
\begin{tikzcd}
\CC \ar[r, shift left, "1"] \ar[r, shift right, "0"'] & \CC \ar[r, shift left, "1"] \ar[r, shift right, "0"'] 
& \CC \ar[r, shift left, "1"] \ar[r, shift right, "0"'] & \CC
\end{tikzcd},\quad
S_-=
\begin{tikzcd}
\CC \ar[r, shift left, "0"] \ar[r, shift right, "1"'] & \CC \ar[r, shift left, "0"] \ar[r, shift right, "1"'] 
& \CC \ar[r, shift left, "0"] \ar[r, shift right, "1"'] & \CC
\end{tikzcd}.
\]

These modules admit the following projective and injective resolutions:
\[
P(4) \overset{\beta_{3}}{\longrightarrow} P(3) \overset{\alpha_{2}}{\longrightarrow} P(2) \overset{\beta_{1}}{\longrightarrow} P(1) \longrightarrow S_+,\quad 
S_+ \longrightarrow I(4) \overset{\tilde{\beta}_{3}}{\longrightarrow} I(3) \overset{\tilde{\alpha}_{2}}{\longrightarrow} I(2) \overset{\tilde{\beta}_{1}}{\longrightarrow} I(1),
\]
\[
P(4) \overset{\alpha_{3}}{\longrightarrow} P(3) \overset{\beta_{2}}{\longrightarrow} P(2) \overset{\alpha_{1}}{\longrightarrow} P(1) \longrightarrow S_-,\quad 
S_- \longrightarrow I(4) \overset{\tilde{\alpha}_{3}}{\longrightarrow} I(3) \overset{\tilde{\beta}_{2}}{\longrightarrow} I(2) \overset{\tilde{\alpha}_{1}}{\longrightarrow} I(1),
\]
where $\tilde{\alpha_i},\tilde{\beta}_i$ denote the homomorphisms $\Hom_{\mod(\CC Q/I)}(\alpha_i,\CC Q/I)^*$, $\Hom_{\mod(\CC Q/I)}(\beta_i,\CC Q/I)^*$.
\end{lem}
\begin{pf}
For simplicity, we only illustrate the case of $S_+$, together with the explicit forms of the maps $\tilde{\alpha_i}$ and $\tilde{\beta_i}$:
\begin{equation*}
\begin{tikzcd}[row sep = 8mm, column sep = 8mm, ampersand replacement=\&]
P(4): \arrow[d, shift right, "\beta_3"] \&[-9mm] 0 \arrow[r, shift left] \arrow[r, shift right] \arrow[d] \& 0 \arrow[r, shift left] \arrow[r, shift right] \arrow[d] \& 0 \arrow[r, shift left] \arrow[r, shift right] \arrow[d] \& \CC \arrow[d, "{\begin{psmallmatrix} 0\\1 \end{psmallmatrix}}"] 
\\
P(3): \arrow[d, shift right, "\alpha_2"] \&[-9mm] 0 \arrow[r, shift left] \arrow[r, shift right] \arrow[d] \& 0 \arrow[r, shift left] \arrow[r, shift right] \arrow[d] \& \CC \arrow[r, shift left] \arrow[r, shift right] \arrow[d, "{\begin{psmallmatrix} 1\\0 \end{psmallmatrix}}"] \& \CC^2 \arrow[d, "{\begin{psmallmatrix} 1&0\\0&0 \end{psmallmatrix}}"] 
\\
P(2): \arrow[d, shift right, "\beta_1"] \&[-9mm] 0 \arrow[r, shift left] \arrow[r, shift right] \arrow[d] \& \CC \arrow[r, shift left] \arrow[r, shift right] \arrow[d, "{\begin{psmallmatrix} 0\\1 \end{psmallmatrix}}"] \& \CC^2 \arrow[r, shift left] \arrow[r, shift right] \arrow[d, "{\begin{psmallmatrix} 0&0\\0&1 \end{psmallmatrix}}"] \& \CC^2 \arrow[d, "{\begin{psmallmatrix} 0&0\\0&1 \end{psmallmatrix}}"] 
\\
P(1): \arrow[d, shift right] \&[-9mm] \CC \arrow[r, shift left] \arrow[r, shift right] \arrow[d, "1"] \& \CC^2 \arrow[r, shift left] \arrow[r, shift right] \arrow[d, "{\begin{psmallmatrix} 1&0 \end{psmallmatrix}}"] \& \CC^2 \arrow[r, shift left] \arrow[r, shift right] \arrow[d, "{\begin{psmallmatrix} 1&0 \end{psmallmatrix}}"] \& \CC^2 \arrow[d, "{\begin{psmallmatrix} 1&0 \end{psmallmatrix}}"] 
\\
S_+: \& \CC \arrow[r, shift left] \arrow[r, shift right] \& \CC \arrow[r, shift left] \arrow[r, shift right] \& \CC \arrow[r, shift left] \arrow[r, shift right] \& \CC
\end{tikzcd}
\quad
\begin{tikzcd}[row sep = 8mm, column sep = 8mm, ampersand replacement=\&]
S_+: \arrow[d, shift right] \&[-9mm] \CC \arrow[r, shift left] \arrow[r, shift right] \arrow[d, "{\begin{psmallmatrix} 1\\0 \end{psmallmatrix}}"] \& \CC \arrow[r, shift left] \arrow[r, shift right] \arrow[d, "{\begin{psmallmatrix} 1\\0 \end{psmallmatrix}}"] \& \CC \arrow[r, shift left] \arrow[r, shift right] \arrow[d, "{\begin{psmallmatrix} 1\\0 \end{psmallmatrix}}"] \& \CC \arrow[d, "1"] 
\\
I(4): \arrow[d, shift right, "\tilde\beta_3"] \&[-9mm] \CC^2 \arrow[r, shift left] \arrow[r, shift right] \arrow[d, "{\begin{psmallmatrix} 0&0\\0&1 \end{psmallmatrix}}"] \& \CC^2 \arrow[r, shift left] \arrow[r, shift right] \arrow[d, "{\begin{psmallmatrix} 0&0\\0&1 \end{psmallmatrix}}"] \& \CC^2 \arrow[r, shift left] \arrow[r, shift right] \arrow[d, "{\begin{psmallmatrix} 0&1 \end{psmallmatrix}}"] \& \CC \arrow[d] 
\\
I(3): \arrow[d, shift right, "\tilde\alpha_2"] \&[-9mm] \CC^2 \arrow[r, shift left] \arrow[r, shift right] \arrow[d, "{\begin{psmallmatrix} 1&0\\0&0 \end{psmallmatrix}}"] \& \CC^2 \arrow[r, shift left] \arrow[r, shift right] \arrow[d, "{\begin{psmallmatrix} 1&0 \end{psmallmatrix}}"] \& \CC \arrow[r, shift left] \arrow[r, shift right] \arrow[d] \& 0 \arrow[d] 
\\
I(2): \arrow[d, shift right, "\tilde\beta_1"] \&[-9mm] \CC^2 \arrow[r, shift left] \arrow[r, shift right] \arrow[d, "{\begin{psmallmatrix} 0&1 \end{psmallmatrix}}"] \& \CC \arrow[r, shift left] \arrow[r, shift right] \arrow[d] \& 0 \arrow[r, shift left] \arrow[r, shift right] \arrow[d] \& 0 \arrow[d] 
\\
I(1): \& \CC \arrow[r, shift left] \arrow[r, shift right] \& 0 \arrow[r, shift left] \arrow[r, shift right] \& 0 \arrow[r, shift left] \arrow[r, shift right] \& 0
\end{tikzcd}
\end{equation*}
The case of $S_-$ is obtained by exchanging every $\alpha_i$ and $\beta_i$.
\qed\end{pf}

\begin{lem}\label{lem:3-sph}
The modules $S_+$ and $S_-$ are $3$-spherical. Moreover, we have 
\[
\RHom_\D(S_+,P(i))\cong \RHom_\D(S_-,P(i))\cong \CC[-3],\quad i=1,\dots,4,
\]
\[
\RHom_\D(S_+,S_i)\cong \RHom_\D(S_-,S_i)\cong \RHom_\D(S_+,S_-)\cong 0, \quad i=1,\dots,3,
\]
where $S_i$'s are the $1$-spherical objects from Proposition~\ref{prop:1-sph}.
\end{lem}
\begin{pf}
This follows from a direct calculation using Lemma~\ref{lem: S_+-} (see also {\cite[Proposition~4.1]{D}}).
\qed\end{pf}

Hence, we obtain the auto-equivalences $T_{S_+}$ and $T_{S_-}$, which commute with $T_{S_{i}}$ for $i=1,\dots, \mu-1$.
It turns out also that $\S_\D[-1]=T_{S_+}\circ T_{S_-}= T_{S_-}\circ T_{S_+}$ in $\Aut(\D)$, but,
since we will not use this fact in this paper, we omit the proof.

In the abelian category ${\rm mod}(\CC Q/I)$, we have short exact sequences
\begin{equation}\label{S+- exact}
0 \longrightarrow P(i) \longrightarrow S_+\oplus S_- \longrightarrow I(i) \longrightarrow 0,\quad i=1,\dots ,4,
\end{equation}
which can be described in terms of quiver representations as follows.
\[
\begin{tikzcd}[row sep = 8mm, column sep = 6mm, ampersand replacement=\&]
P(1): \arrow[d, shift right] \&[-9mm] \CC \arrow[r, shift left] \arrow[r, shift right] \arrow[d, "{\begin{psmallmatrix} 1\\-1 \end{psmallmatrix}}"] \& \CC^2 \arrow[r, shift left] \arrow[r, shift right] \arrow[d, "{\begin{psmallmatrix} 1&0\\0&-1 \end{psmallmatrix}}"] \& \CC^2 \arrow[r, shift left] \arrow[r, shift right] \arrow[d, "{\begin{psmallmatrix} 1&0\\0&-1 \end{psmallmatrix}}"] \& \CC^2 \arrow[d, "{\scriptstyle \begin{psmallmatrix} 1&0\\0&-1 \end{psmallmatrix}}"] 
\\
S_+\oplus S_-: \arrow[d, shift right] \& \CC^2 \arrow[r, shift left] \arrow[r, shift right] \arrow[d, "{\begin{psmallmatrix} 1&1 \end{psmallmatrix}}"] \& \CC^2 \arrow[r, shift left] \arrow[r, shift right] \arrow[d] \& \CC^2 \arrow[r, shift left] \arrow[r, shift right] \arrow[d] \& \CC^2 \arrow[d] 
\\
I(1): \& \CC \arrow[r, shift left] \arrow[r, shift right] \& 0 \arrow[r, shift left] \arrow[r, shift right] \& 0 \arrow[r, shift left] \arrow[r, shift right] \& 0
\end{tikzcd}
\
\begin{tikzcd}[row sep = 8mm, column sep = 6mm, ampersand replacement=\&]
P(2): \arrow[d, shift right] \&[-9mm] 0 \arrow[r, shift left] \arrow[r, shift right] \arrow[d] \& \CC \arrow[r, shift left] \arrow[r, shift right] \arrow[d, "{\begin{psmallmatrix} 1\\-1 \end{psmallmatrix}}"] \& \CC^2 \arrow[r, shift left] \arrow[r, shift right] \arrow[d, "{\begin{psmallmatrix} 1&0\\0&-1 \end{psmallmatrix}}"] \& \CC^2 \arrow[d, "{\begin{psmallmatrix} 1&0\\0&-1 \end{psmallmatrix}}"] 
\\
S_+\oplus S_-: \arrow[d, shift right] \& \CC^2 \arrow[r, shift left] \arrow[r, shift right] \arrow[d, "{\begin{psmallmatrix} 1&0\\0&1 \end{psmallmatrix}}"] \& \CC^2 \arrow[r, shift left] \arrow[r, shift right] \arrow[d, "{\begin{psmallmatrix} 1&1 \end{psmallmatrix}}"] \& \CC^2 \arrow[r, shift left] \arrow[r, shift right] \arrow[d] \& \CC^2 \arrow[d] 
\\
I(2): \& \CC^2 \arrow[r, shift left] \arrow[r, shift right] \& \CC \arrow[r, shift left] \arrow[r, shift right] \& 0 \arrow[r, shift left] \arrow[r, shift right] \& 0
\end{tikzcd}
\]
\[
\begin{tikzcd}[row sep = 8mm, column sep = 6mm, ampersand replacement=\&]
P(3): \arrow[d, shift right] \&[-9mm] 0 \arrow[r, shift left] \arrow[r, shift right] \arrow[d] \& 0 \arrow[r, shift left] \arrow[r, shift right] \arrow[d] \& \CC \arrow[r, shift left] \arrow[r, shift right] \arrow[d, "{\begin{psmallmatrix} 1\\-1 \end{psmallmatrix}}"] \& \CC^2 \arrow[d, "{\begin{psmallmatrix} 1&0\\0&-1 \end{psmallmatrix}}"] 
\\
S_+\oplus S_-: \arrow[d, shift right] \& \CC^2 \arrow[r, shift left] \arrow[r, shift right] \arrow[d, "{\begin{psmallmatrix} 1&0\\0&1 \end{psmallmatrix}}"] \& \CC^2 \arrow[r, shift left] \arrow[r, shift right] \arrow[d, "{\begin{psmallmatrix} 1&0\\0&1 \end{psmallmatrix}}"] \& \CC^2 \arrow[r, shift left] \arrow[r, shift right] \arrow[d, "{\begin{psmallmatrix} 1&1 \end{psmallmatrix}}"] \& \CC^2 \arrow[d] 
\\
I(3): \& \CC^2 \arrow[r, shift left] \arrow[r, shift right] \& \CC^2 \arrow[r, shift left] \arrow[r, shift right] \& \CC \arrow[r, shift left] \arrow[r, shift right] \& 0
\end{tikzcd}
\
\begin{tikzcd}[row sep = 8mm, column sep = 6mm, ampersand replacement=\&]
P(4): \arrow[d, shift right] \&[-9mm] 0 \arrow[r, shift left] \arrow[r, shift right] \arrow[d] \& 0 \arrow[r, shift left] \arrow[r, shift right] \arrow[d] \& 0 \arrow[r, shift left] \arrow[r, shift right] \arrow[d] \& \CC \arrow[d, "{\begin{psmallmatrix} 1\\-1 \end{psmallmatrix}}"] 
\\
S_+\oplus S_-: \arrow[d, shift right] \& \CC^2 \arrow[r, shift left] \arrow[r, shift right] \arrow[d, "{\begin{psmallmatrix} 1&0\\0&1 \end{psmallmatrix}}"] \& \CC^2 \arrow[r, shift left] \arrow[r, shift right] \arrow[d, "{\begin{psmallmatrix} 1&0\\0&1 \end{psmallmatrix}}"] \& \CC^2 \arrow[r, shift left] \arrow[r, shift right] \arrow[d, "{\begin{psmallmatrix} 1&0\\0&1 \end{psmallmatrix}}"] \& \CC^2 \arrow[d, "{\begin{psmallmatrix} 1&1 \end{psmallmatrix}}"] 
\\
I(4): \& \CC^2 \arrow[r, shift left] \arrow[r, shift right] \& \CC^2 \arrow[r, shift left] \arrow[r, shift right] \& \CC^2 \arrow[r, shift left] \arrow[r, shift right] \& \CC
\end{tikzcd}
\]
Suppose that $\E_P$ and $T_{S_+}^k\!\cdot \E_P$ are in the same $B_4\ltimes \ZZ^4$-orbit for some non-zero integer $k$.
Then, by Proposition~\ref{prop:l-r braid}, there would exist $\Phi\in\langle T_{S_1},\dots,T_{S_4}\rangle$
and $n\in\ZZ$ such that $(\Phi\circ[n])\cdot \E_P\cong T_{S_+}^k\cdot\E_P$.
Namely, we have $(T_{S_+}^{-k}\circ \Phi\circ [n])(P(i))\cong P(i)$ for $i=1,\dots, 4$. 
Hence the auto-equivalence $T_{S_+}^{-k}\circ \Phi\circ [n]$ respects the standard $t$-structure of $\D$ and we have $(T_{S_+}^{-k}\circ \Phi\circ [n])(\mod(\CC Q/I))=\mod(\CC Q/I)$.
Recall that the Serre functor $\S_\D$ commutes with any auto-equivalence of $\D$ and that $I(i)\cong\S_\D(P(i))$, $i=1,\dots, 4$.
However, we have
\[
(T_{S_+}^{-k}\circ \Phi\circ [n])(S_+\oplus S_-)\cong T_{S_+}^{-k}(\Phi(S_+[n]\oplus S_-[n]))\cong T_{S_+}^{-k}(S_+[n]\oplus S_-[n])\cong S_+[n-2k]\oplus S_-[n],
\]
which contradicts \eqref{S+- exact} that $S_+, S_-\in \mod(\CC Q/I)$ since $k\ne 0$. 
In particular, it is concluded that $\FEC(\D)$ has an infinitely many $B_4\ltimes \ZZ^4$-orbits (at least, it has orbits indexed by $k\in \ZZ$).

\subsection{The case $\mu\geq 5$}

We next prove Theorem~\ref{thm: main} for the case $\mu\geq 5$.

Consider the full strong exceptional collection $\E_P=(P(\mu),\dots,P(1))$ in $\D$ and set $P_4:=P(4)\oplus\dots \oplus P(1)$.
Since the $\CC$-algebra $\End_\D(P_4):=\Hom_\D(P_4,P_4)$ is isomorphic to the path algebra $\CC Q/I$ with $\mu=4$, the functor 
\[
\RHom_\D(P_4,-):\D\longrightarrow \D^{b}\mod(\End_\D(P_4)),\quad X\mapsto \RHom_\D(P_4,X),
\]
gives an equivalence between the smallest full triangulated subcategory $\langle P(4),\dots,P(1)\rangle$ of $\D$ containing $P(4),\dots,P(1)$ 
and the derived category $\D^{b}\mod(\CC Q/I)$ with $\mu=4$.
In particular,  under this equivalence, ``$P(4),\dots,P(1)$ in $\D$'' and ``$P(4),\dots,P(1)$ in $\D^{b}\mod(\CC Q/I)$ with $\mu=4$'' correspond literally and can be identified.
Therefore, for $k\in\ZZ$, through this identification a full exceptional collection $(T_{S_+}^k(P(4)),\dots,T_{S_+}^k(P(1)))$ in $\D^{b}\mod(\CC Q/I)$ with $\mu=4$  
can be completed to a full exceptional collection $\E'_k$ in $\D$ as 
\begin{equation*}
\E'_k := (P(\mu),\dots,P(5),T_{S_+}^k(P(4)),\dots,T_{S_+}^k(P(1))),\quad \E'_0=\E_P,
\end{equation*}
where the same letters represent the corresponding objects in $\D$.
For each $i=1,\dots,4$, we have the exact triangle
\begin{equation}\label{3-ST exact tri}
S_+[-2k-3] \longrightarrow T_{S_+}^{k}(P(i)) \longrightarrow T_{S_+}^{k+1}(P(i)) \longrightarrow S_+[-2k-2]
\end{equation}
in $\D$ by the definition of the functor $T_{S_+}$ on $\D^{b}\mod(\CC Q/I)$ with $\mu=4$, 
again through the above identification.
\begin{lem}
For $i=1,\dots,4$, $j=5,\dots,\mu$ and $k\in\ZZ$, we have
\[
\RHom_\D\big(P(j),\,T_{S_+}^k(P(i))\big)\cong\CC\oplus\CC[-2k].
\]
\end{lem}
\begin{pf}
If $k=0$, then the statement is obvious. 
Lemma~\ref{lem: S_+-} yields that $S_+$ is ismorphic in $\D$ to the complex $P(4) \overset{\beta_{3}}{\longrightarrow} P(3) \overset{\alpha_{2}}{\longrightarrow} P(2) \overset{\beta_{1}}{\longrightarrow} P(1)$ where the term $P(1)$ has degree zero.
Therefore, it follows that $\RHom_\D(P(j),S_+)\cong \CC\oplus\CC[3]$.
Applying $\RHom_\D(P(j),-)$ to the exact triangle \eqref{3-ST exact tri}, 
we obtain the statement inductively.
\qed\end{pf}
Hence we see that $\E'_k$ is \emph{not} strong if $k\ne 0$. 
Moreover, $\E'_k$ cannot be transferred to a full strong exceptional collection by the action of $\ZZ^\mu$.
On the other hand, $\E_P$ is strong, and the action of $B_\mu$ respects the property of being strong by Proposition~\ref{prop:l-r braid}.
Consequently, if $k\ne 0$ then $\E'_k$ cannot be reachable from $\E_P$ under the $B_\mu\ltimes \ZZ^\mu$-action.
This completes the proof of Theorem~\ref{thm: main}.


\end{document}